\documentclass[12pt,reqno]{amsart}
\usepackage{enumerate, latexsym, amsmath, amsfonts, amssymb, amsthm, color}
\def\pmod #1{\ ({\rm{mod}}\ #1)}
\def\Z{\Bbb Z}
\def\N{\Bbb N}

\def\Q{\Bbb Q}

\def\l{\left}
\def\r{\right}
\def\bg{\bigg}
\def\({\bg(}
\def\){\bg)}
\def\t{\text}
\def\f{\frac}

\def\ls{\leqslant}
\def\gs{\geqslant}

\def\sm{\setminus}

\def\al{\alpha}

\def\eq{\equiv}

\def\Proof{\noindent{\it Proof}}

\theoremstyle{plain}
\newtheorem{theorem}{Theorem}

\newtheorem{lemma}{Lemma}

\theoremstyle{definition}

\theoremstyle{remark}
\newtheorem{remark}{Remark}

 \vspace{4mm}

\begin{document}

\hbox{Preprint}
\medskip

\title
[{Diophantine equations over $\mathbb Z[i]$ with 20 unknowns}]
{Undecidability on Diophantine equations\\ over $\mathbb Z[i]$ with $20$ unknowns}

\author
[Y. Matiyasevich and Z.-W. Sun] {Yuri Matiyasevich and Zhi-Wei Sun}

\address {(Yuri Matiyasevich) St. Petersburg Department of Steklov Mathematical Institute of Russian Academy of Sciences, Fontanka 27, 191023, St. Petersburg, Russia}
\email{yumat@pdmi.ras.ru}

\address{(Zhi-Wei Sun, corresponding author) School of Mathematics, Nanjing
University, Nanjing 210093, People's Republic of China}
\email{zwsun@nju.edu.cn}

\subjclass[2010]{Primary 11U05, 03D35; Secondary 03D25, 11D99, 11R11.}
\keywords{Hilbert's Tenth Problem, Diophantine equation, Gaussian ring, imaginary quadratic field, undecidability.
\newline \indent The second author is supported by the Natural Science Foundation of China (grant no. 12371004).}

\begin{abstract} It is known that Hilbert's Tenth Problem over the Gaussian ring $\mathbb Z[i]=\{a+bi:\ a,b\in\Z\}$ is undecidable. In this paper we obtain the following further result: There is no algorithm to decide whether
an arbitrarily given polynomial equation $P(z_1,\ldots,z_{20})$ $=0$
(with integer coefficients) is solvable over $\mathbb Z[i]$.
This improves the previous record involving $52$ variables.
\end{abstract}
\maketitle

\section{Introduction}
\setcounter{lemma}{0}
\setcounter{theorem}{0}
\setcounter{corollary}{0}
\setcounter{remark}{0}

Hilbert's Tenth Problem (HTP in short) over a ring $R$ asks for an algorithm  to test whether an arbitrary polynomial Diophantine equation
$$P(x_1,\ldots,x_n)=0$$
with coefficients in $R$ has solutions over $R$.
This was posed in 1900 by D. Hilbert when $R$ is the ring $\Z$ of integers.
In 1961, M. Davis, H. Putnam and J. Robinson
\cite{DPR} proved that the solvability of exponential Diophantine equations over $\N$
is undecidable. The original HTP over $\Z$ was finally solved by Y. Matiyasevich \cite{M70} negatively in 1970, see also the book \cite{M93}. In this direction, the second author \cite{S17} proved further that there is no algorithm
to decide for any given $P(x_1,\ldots,x_{11})\in\Z[x_1,\ldots,x_{11}]$ whether the equation
$P(x_1,\ldots,x_{11})=0$ has integer solutions.

Let $K$ be any number field (which is a finite extension of the field $\Q$ of rational numbers),
and let $O_K$ be the ring of algebraic integers in $K$.
What about HTP over $O_K$? If $\Z$ is Diophantine over $O_K$, i.e., there is a polynomial
$P(x_0,x_1,\ldots,x_n)$ over $O_K$ such that $a\in O_K$ lies in $\Z$ if and only if
$P(a,x_1,\ldots,x_n)=0$ for some $x_1,\ldots,x_n\in O_K$, then
HTP over $O_K$ is undecidable via  Matiyasevich's theorem.
J. Denef \cite{D75,D80} proved that if $K$ is an imaginary quadratic field or
a totally real field then $\Z$ is Diophantine over $O_K$ and hence HTP over $O_K$ is unsolvable
(see also \cite[pp.\,98-100]{Sh}).
Recently, P. Koymans and C. Pagano \cite{KP}, as well as L. Alp\"oge, M. Bhargava, W. Ho
and A. Shnidman \cite{ABHS}
successfully proved that $O_K$ is always Diophantine over $\Z$ and hence
HTP over $O_K$ is unsolvable.

The authors \cite{MS} proved that there is no algorithm to decide whether
for any $P(x_1,\ldots,x_{52})\in\Z[x_1,\ldots,x_{52}]$ the  equation
$P(x_1,\ldots,x_{52})=0$ is solvable over $\Z[i]$.
In this paper we make further progress by establishing the following result.

\begin{theorem} \label{Th1.1} There is no algorithm to decide for any $P(x_1,\ldots,x_{20})\in\Z[x_1,\ldots,x_{20}]$
whether the diophantine equation
$$P(x_1,\ldots,x_{20})=0$$
with $20$ unknowns has solutions with $x_1,\ldots,x_{20}\in\Z[i]$.
\end{theorem}

The following auxiliary theorem plays an important role in our proof of Theorem \ref{Th1.1}.

\begin{theorem} \label{Th1.2} Let $d$ be a squarefree positive integer. Let $K$ be the imaginary quadratic field
$\Q(\sqrt {-d})$, and let $O_K$ be the ring of algebraic integers in $K$. Let $x_1,\ldots,x_n\in O_K$
and set $y=2\prod_{k=1}^n(3x_k+1)$. Then we have
$$y+\sum_{k=1}^n\f{x_k}{y^k}\in\Q\iff x_1,\ldots,x_n\in\Z.$$
\end{theorem}

We are going to prove Theorem \ref{Th1.2} and Theorem \ref{Th1.1} in Sections 2 and 3, respectively.

In 1934 T. Skolem \cite{Sk} used a simple trick to reduce any polynomial diophantine equation over $\Z$ to one with degree $4$. By the same trick and Lemma \ref{Lem3.2}, we can reduce any polynomial diophantine equation over $\Z[i]$ to one with degree $4$. Thus the solvability of a general polynomial
equation of degree four over $\Z[i]$ is undecidable.

In 1972 C. L. Siegel \cite{Si} proved that
the solvability of a general quadratic equation (with integer coefficients) over $\Z$ is decidable!
Now it is natural to ask the following question which might have a negative answer.
\medskip

\noindent {\bf Question 1.1}. {\it Whether the solvability of a general quadratic equation over $\Z[i]$ is decidable?}

\section{Proof of Theorem \ref{Th1.2}}
\setcounter{lemma}{0}
\setcounter{theorem}{0}
\setcounter{corollary}{0}
\setcounter{remark}{0}
\setcounter{equation}{0}

Let $d$ be a squarefree positive integer, and set $K=\Q(\sqrt {-d})$.
For $\al\in K$, we write $\al=r(\al)+s(\al)\sqrt {-d}$ with $r(\al),s(\al)\in\Q$, and the norm of $\al$
is given by $N(\al)=r(\al)^2+d\,s(\al)^2$. It is well known that $N(\al\beta)=N(\al)N(\beta)$ for all $\al,\beta\in K$.
Note that
$$|s(\al)|\ls\f{N(\al)}{\sqrt{d}}\qquad \ \t{for any}\ \al\in K.$$
 It is well known that
$$O_K=\begin{cases}\{a+b\sqrt {-d}:\ a,b\in\Z\}&\t{if}\ d\not\eq-1\pmod4,
\\\{\f{a+b\sqrt {-d}}2:\ a,b\in\Z\ \t{and}\ 2\mid a-b\}&\t{if}\ d\eq-1\pmod4.
\end{cases}$$
For $x\in O_K$, we clearly have $N(x)\in\N=\{0,1,2,\ldots\}$, and $N(x)=0$ if and only if $x=0$.

\begin{lemma}\label{Lem2.1} For any $x\in O_K$ we have $N(3x+1)\gs N(x).$
\end{lemma}
\Proof. Write $x=(a+b\sqrt {-d})/2$ with $a,b\in\Z$ and $a\eq b\pmod2$. Then
$$N(3x+1)=N\l(\f{3a+2+3b\sqrt {-d}}2\r)=\f{(3a+2)^2+d(9b^2)}4\gs\f{(3a+2)^2+db^2}4.$$
So it suffices to prove the inequality $(3a+2)^2\gs a^2$. Note that
$$(3a+2)^2-a^2=4(a+1)(2a+1)\gs0$$
since $a\not\in(-1/2,-1)$. This ends our proof. \qed

\begin{lemma}\label{Lem2.2} Let $x_1,\ldots,x_n\in O_K$ and $y=2x_0\prod_{k=1}^n(3x_k+1)$
with $x_0\in O_K\sm\{0\}$. Let $z\in O_K$ with $z+\sum_{k=1}^n x_k/y^k\in\Q$. Then $z\in\Z$.
\end{lemma}
\Proof. Let $Z=\sum_{k=1}^n x_k/y^k$. As $z+Z\in\Q$, we have  $s(z)=-s(Z).$ Note that $N(x_0)\in\Z^+=\{1,2,3,\ldots\}$.
For each $k=1,\ldots,n$, clearly $N(3x_k+1)\in\Z^+$
and
$$1\ls N(3x_k+1)\ls N(x_0)\prod_{j=1}^n N(3x_j+1)=N\l(\f y2\r)=\f{N(y)}4,$$
hence $N(x_k)\ls N(y)/4$ by Lemma \ref{Lem2.1}.
Observe that
\begin{equation}\label{|s(z)|}|s(z)|=|s(Z)|\ls \sum_{k=1}^n \l|s\l(\f{x_k}{y^k}\r)\r|\ls\sum_{k=1}^n\sqrt{\f{N(x_k/y^k)}{d}}
=\f1{\sqrt{d}}\sum_{k=1}^n\sqrt{\f{N(x_k)}{N(y)^k}}
\end{equation}
and thus
\begin{equation}\label{s(y)<}\sqrt{d}\,|s(z)|\ls\sum_{k=1}^n\f1{2\sqrt{N(y)}^{k-1}}\ls\sum_{k=1}^n\f1{2^k}
<\sum_{k=1}^\infty\f1{2^k}=1.
\end{equation}

Note that $z\in O_K$ and $2s(z)\in\Z$. We want to show $z\in\Z=\Q\cap O_K$ (i.e., $s(z)=0$).
If $d\gs5$, then by \eqref{s(y)<} we have
$$2|s(z)|\ls\sqrt d\,|s(z)|<1$$
and hence $s(z)=0$.
For $d\in\{1,2\}$, by \eqref{s(y)<} we have
$$|s(z)|\ls\sqrt d\,|s(z)|<1$$
and hence $s(z)=0$ since $s(z)\in\Z$.

Now we handle the case $d=3$. Without loss of generality, we may simply assume that $x_1,\ldots,x_n$ are nonzero. Let $\omega$ denote the cubic root $(-1+\sqrt{-3})/2$ of unity.
Then $\pm1,\pm\omega,\pm\omega^2$
are the only units in $O_K=\{a+b\omega:\ a,b\in\Z\}.$
 For $\al=a+b\omega\not=0$ with $a,b\in\Z$, clearly
$3\al+1=3a+1+3b\omega\not\in\{\pm1,\pm\omega,\pm\omega^2\}$ and hence $N(3\al+1)\not=1$. Since $x_1,\ldots,x_n\in O_K\sm\{0\}$, we have $N(3x_k+1)\gs2$ for all $k=1,\ldots,n$.
Thus, for each $1\ls k\ls n$, we have
$$N\l(\f y2\r)=N(x_0) N(3x_k+1)\prod_{j=1\atop j\not=k}^nN(3x_j+1)\gs 2^{k-1}N(3x_k+1)\gs 2^{k-1}N(x_k)$$
and hence $N(y)\gs 2^{k+1}N(x_k)\gs2^{k+1}$.
Combining this with \eqref{|s(z)|}, we obtain that
\begin{align*}\sqrt3\,|s(z)|&\ls \sum_{k=1}^n\sqrt{\f{N(x_k)}{N(y)^k}}
\ls\sum_{k=1}^n\sqrt{\f{N(y)/2^{k+1}}{N(y)^k}}=\f12\sum_{k=1}^n\f1{\sqrt{2N(y)}^{k-1}}
\\&\ls\f12\sum_{k=1}^n\f1{2^{(k+2)(k-1)/2}}<\f12\l(1+\f1{2^2}+\f1{2^5}+\sum_{k=4}^\infty\f1{2^{k-1}}\r)
=\f{49}{64}
\end{align*}
and hence
$$|s(z)|<\f{49}{64\sqrt3}<\f12.$$
Thus $s(z)=0$ since $2s(z)\in\Z$.

In view of the above, we have completed the proof of Lemma \ref{Lem2.2}. \qed

\medskip
\noindent{\tt Proof of Theorem \ref{Th1.2}}. Set $Y=\sum_{k=1}^n x_k/y^k$.

 If $x_1,\ldots,x_n\in \Z$, then
$y=2\prod_{k=1}^n(3x_k+1)\in\Z$ and hence $y+Y\in\Q$.

Now suppose $y+Y\in\Q$. By Lemma \ref{Lem2.2}, $y\in\Z$ and thus $Y\in\Q$.
We want to use induction to show that $x_m\in\Z$ for all $m=1,\ldots,n$.
Fix $m\in\{1,\ldots,n\}$, and assume that $x_k\in\Z$ for all $0<k<m$.
Then
$$\sum_{k=m}^n\f {x_k}{y^k}=Y-\sum_{0<k<m}\f{x_k}{y^k}\in\Q$$
and hence
\begin{equation}\label{xm}x_m+\sum_{m<k\ls n}\f{x_k}{y^{k-m}}=y^m\sum_{k=m}^n\f {x_k}{y^k}\in\Q.
\end{equation}
Observe that $y=2x_0\prod_{m<k\ls n}(3x_k+1)$ with $x_0=\prod_{1\ls k\ls m}(3x_k+1)\in O_K\sm\{0\}$.
In view of \eqref{xm}, by applying Lemma \ref{Lem2.2} we obtain $x_m\in\Z$.
This concludes the induction step.

By the above, we have completed the proof of Theorem \ref{Th1.2}. \qed

\section{Proof of Theorem \ref{Th1.1}}
\setcounter{lemma}{0}
\setcounter{theorem}{0}
\setcounter{corollary}{0}
\setcounter{remark}{0}
\setcounter{equation}{0}

Motivated by the Matiyasevich--Robinson Relation-Combining Theorem (cf. \cite{MR}), the second author
 \cite[p.\,69]{S-book} deduced the following lemma for the rational field $\Q$.
 We now extend it to any number field.

\begin{lemma}\label{Lem3.1} Let $K$ be a number field, and let $A_1,A_2,S,T\in O_K$ with $A_1\not=A_2$
and $S\not=0$. Then
$$A_1\in\square\land A_2\in\square\land S\mid T$$
if and only if for some $m\in O_K$ we have
\begin{equation}\label{ST} f(A_1,A_2,S,T,m)=0
\end{equation}
where $\square=\{\al^2:\ \al\in O_K\}$ and
$$f(A_1,A_2,S,T,m):=(T-mS)^4-2(A_1+A_2)S^2(T-mS)^2+(A_1-A_2)^2S^4.$$
\end{lemma}
\Proof. Observe that
\begin{align*}&\ S^{-4}f(A_1,A_2,S,T,m)
\\=&\ \l(\l(\f TS-m\r)^2-(A_1+A_2)\r)^2-4A_1A_2
\\=&\ \l(\l(\f TS-m\r)^2+A_1-A_2\r)^2-4A_1\l(\f TS-m\r)^2.
\end{align*}

If $A_1=\al_1^2$ and $A_2=\al_2^2$ with $\al_1,\al_2\in O_K$, and $S\mid T$,
then $m=T/S-\al_1-\al_2\in O_K$ and
\begin{align*}&\ \l(\l(\f TS-m\r)^2-(A_1+A_2)\r)^2
\\=&\ ((\al_1+\al_2)^2-\al_1^2-\al_2^2)^2=(2\al_1\al_2)^2=4A_1A_2,
\end{align*}
thus $f(A_1,A_2,S,T,m)=0$.

Now assume that \eqref{ST} holds for some $m\in O_K$. Set $x=T/S-m$. As $(x^2+A_1-A_2)^2-4A_1x^2=0$,
$x$ is an algebraic integer and hence $x\in O_K$. Thus $S\mid T$. Since $A_1\not=A_2$, we have $x\not=0$. For $y=(A_1-A_2)/x\in K$, we have $(x+y)^2=4A_1$ and hence $y$ is an algebraic integer.
As $x,y\in O_K$ and $(\f{x+y}2)^2=A_1\in O_K$, $\f{x+y}2$ must be an algebraic integer.
Note that
$$A_2=A_1-xy=\l(\f{x+y}2\r)^2-xy=\l(\f{x-y}2\r)^2=\l(\f{x+y}2-y\r)^2.$$
So $A_1,A_2\in\square$.

By the above, we have finished the proof of Lemma \ref{Lem3.1}. \qed

The following lemma is an easy fact, see, e.g., \cite[Lemma 2.4]{MS}.

\begin{lemma} \label{Lem3.2} \label{land} For any $x,y\in\Z[i]$, we have
$$x=0 \land y=0\iff x^2+2y^2=0.$$
\end{lemma}

\begin{lemma}\label{Lem3.3} A number $t\in\Z[i]$ is a rational integer if and only if
there are $v,x,y\in\Z[i]$ with $v\not=0$ such that
\begin{equation}\label{zzi}
4(2v(2(2t+1)^2+1)-y)^2-3y^2-1=0
\end{equation}
and
\begin{equation} 3y^2(2t+1-xy)^2+1\in\square,
\end{equation}
where $\square=\{\al^2:\ \al\in\Z[i]\}$. When $t\in\Z$, we can actually require further that $v,x,y\in\Z$.
\end{lemma}
\begin{remark} Lemma \ref{Lem3.3} follows from  \cite[Theorem 1.1]{MS} and its proof in view of Lemma \ref{land}.
\end{remark}

We also need the following result observed by S. P. Tung \cite{T85},

\begin{lemma} \label{Lem3.4} An integer $m$ is nonzero if and only if $m=(2r+1)(3s+1)$ for some $r,s\in\Z$.
\end{lemma}

\medskip
\noindent{\tt Proof of Theorem \ref{Th1.1}}. Let  ${\mathcal A}$ be a subset of $\N$ which is recursively enumerable but not recursive. (The existence of such a set is well known, see, e.g., N. Cutland \cite[pp.\,140-141]{C80}.) By Sun \cite[Theorem 1.1(ii)]{S17}, there is a
polynomial $P(z_0,\ldots,z_{10})\in\Z[z_0,\ldots,z_{10}]$ such that $a\in\N$ belongs to ${\mathcal A}$
if and only if
\begin{equation}\label{P} P(a,z_1,\ldots,z_{10})=0
 \end{equation} for some $z_1,\ldots,z_{10}\in\Z$ with $z_{10}\not=0$.

Let $y=2\prod_{k=1}^{10}(3z_k+1)$. By Theorem \ref{Th1.1}, when $z_1,\ldots,z_{10}\in\Z[i]$, we have
\begin{align*}&\ z_1,\ldots,z_{10}\in\Z
\\\iff&\ y+\sum_{k=1}^{10}\f{z_k}{y^k}\in \Q
\\\iff&\ \exists t\in\Z\ \l[t\not=0\land t\l(y+\sum_{k=1}^{10}\f{z_k}{y^k}\r)\in\Z\r]
\\\iff&\ \exists t\in\Z\ \l[t\not=0\land ty^{10}\mid \sum_{k=1}^{10}z_ky^{10-k}\land t\l(y+\sum_{k=1}^{10}\f{z_k}{y^k}\r)\in\Z\r].
\end{align*}

Suppose that $t\in\Z[i]$ and $ty^{10}\mid \sum_{k=1}^{10}z_ky^{10-k}$.
By Lemma \ref{Lem3.3}, $t\in\Z$ if and only there are $v_0,x_0,y_0\in\Z[i]$ with $v_0\not=0$ such that
\begin{equation}\label{v0}4(2v_0(2(2t+1)^2+1)-y_0)^2-3y_0^2-1=0
\end{equation}
and
\begin{equation} A_1:=9(3y_0^2(2t+1-x_0y_0)^2+1)\in\square,
\end{equation}
and
$t(y+\sum_{k=1}^{10}z_k/{y^k})\in\Z$ if and only if there are $v_1,x_1,y_1\in\Z[i]$ with $v_1\not=0$ such that
\begin{equation}\label{v1}4(2v_1(2(2t(y^{11}+z_1y^9+\cdots+z_{10}y^0)+y^{10})^2+y^{20})-y_1y^{20})^2
=(3y_1^2+1)y^{40}
\end{equation}
and
\begin{equation}A_2:=3y_1^2(2t(y^{11}+z_1y^9+\ldots+z_{10}y^0)+y^{10}-x_1y_1y^{10})^2+y^{20}\in\square.
\end{equation}
Moreover, when $z_1,\ldots,z_{10},t\in\Z$, we can require further that $v_0,x_0,y_0\in\Z$,
and also that $v_1,x_1,y_1\in\Z$ if $t(y+\sum_{k=1}^{10}z_k/{y^k})\in\Z$.
Note that $A_1\not=A_2$ since $A_1\eq0\pmod3$ and $A_2\eq y^{20}\eq2^{20}\eq1\pmod{3}$.

Let $a\in\N$. By the above, $a\in \mathcal A$ if and only if there are $$z_1,\ldots,z_{10},t,v_0,x_0,y_0,v_1,x_1,y_1\in\Z[i]$$
with $z_{10}tv_0v_1\not=0$ such that \eqref{P}, \eqref{v0} and \eqref{v1} all hold and also
$$A_1\in\square \land A_2\in\square \land ty^{10}\mid \sum_{k=1}^{10}z_ky^{10-k}.$$
Note that when $a\in A$ we can actually find
$$z_1,\ldots,z_{10},t,v_0,x_0,y_0,v_1,x_1,y_1\in\Z$$
to meet the requirement. When $z_{10}tv_0v_1\in\Z\sm\{0\}$, there are $r,s\in\Z$
such that $z_{10}tv_0v_1=(2r+1)(3s+1)$. If the last equality holds for some $r,s\in\Z[i]$,
then we obviously have $z_{10}tv_0v_1\not=0$. Thus, in view of Lemma \ref{Lem3.1},
$a\in \mathcal A$ if and only if there are
$$z_1,\ldots,z_{10},m,r,s,t,v_0,x_0,y_0,v_1,x_1,y_1\in\Z[i]$$
such that \eqref{P}, \eqref{v0}, \eqref{v1}, and the equalities
$$z_{10}tv_0v_1=(2r+1)(3s+1)$$
and
$$f\l(A_1,A_2,ty^{10}, \sum_{k=1}^{10}z_ky^{10-k},m\r)=0$$
all hold. Therefore, by using Lemma \ref{Lem3.2} we see that
$a\in\mathcal A$ if and only if
$$F(a,z_1,\ldots,z_{10},m,r,s,t,v_0,x_0,y_0,v_1,x_1,y_1)=0$$
for some
$$z_1,\ldots,z_{10},m,r,s,t,v_0,x_0,y_0,v_1,x_1,y_1\in\Z[i],$$
where $F$ is a suitable polynomial with integer coefficients in 21 variables.

As $\mathcal A$ is not recursive, by the last paragraph we have the desired result.
\qed


\begin{thebibliography}{DPR}

\bibitem{ABHS} L. Alp\"oge, M. Bhargava, W. Ho
and A. Shnidman, {\it Rank stability in quadratic extensions and Hilbert's tenth problem
for the ring of integers of a number field}, arXiv:2501.18774, 2025.

\bibitem{C80} N. Cutland,  Computability, Cambridge Univ. Press, Cambridge, 1980.

\bibitem{DPR} M. Davis, H. Putnam and J. Robinson, {\it The decision problem for exponential diophantine equations}, Ann. of Math. {\bf 74}(1961), 425--436.

\bibitem{D75} J. Denef, {\it Hilbert's Tenth Problem for quadratic rings}, Proc. Amer. Math. Soc.
{\bf 48} (1975), 214--220.

\bibitem{D80} J. Denef, {\it Diophantine sets of algebraic integers, II},
Trans. Amer. Math. Soc. {\bf 257 (1980)}, 227--236, 1980.

\bibitem{KP} P. Koymans and C. Pagano, Hilbert's tenth problem via additive combinatorics,
arXiv:2412.01768, 2024.

\bibitem{M70} Y.  Matiyasevich, {\it Enumerable sets are diophantine}, Dokl. Akad. Nauk SSSR {\bf 191} (1970), 279--282;
English translation with addendum, Soviet Math. Doklady {\bf 11} (1970), 354--357.

\bibitem{M93} Y. Matiyasevich, Hilbert's Tenth Problem, MIT Press, Cambridge, Massachusetts, 1993.
Translation from Russian original, Nauka Publisher, Moscow.

\bibitem {MR} Y Matiyasevich and J. Robinson, {\it Reduction of an arbitrary diophantine equation to one in 13 unknowns},  Acta Arith. {\bf 27} (1975), 521--553.

\bibitem{MS} Y. Matiyasevich and  Z.-W. Sun, {\it On Diophantine equations over $\mathbb Z[i]$ with $52$ unknowns},
in: Mathematical Logic, Computability, Complexity and Randomness
(edited by J. Brendle et al.), pp. 153--158, World Sci., 2025.

\bibitem{Si} C. L. Siegel, {\it Zur Theorie der quadratischen Formen}, Nachrichten der Akademie der Wissenschaften in G\"ottingen II. Mathematisch-Physikalische Klasse, {\bf 3} (1934), 21--46.

\bibitem{Sh} A. Shlapentokh, Hilbert's Tenth Problem: Diophantine Classes and Extensions to Global Fields, New Mathematical Monographs, Vol. 7, Cambridge Univ. Press, Cambridge, 2007.
    
\bibitem{Sk} T. Skolem, {\it \"Uber die Nicht-charakterisierbarkeit der Zahlenreihe mittels endlich oder abz\"ahlbar unendlich vieler Aussagen mit ausschliesslich Zahlenvariablen},
    Fund. Math. {\bf 23} (1934), 150--161.

\bibitem{S17} Z.-W. Sun, {\it Further results on Hilbert's tenth problem}, Sci. China Math.
{\bf 64} (2021), 281--306.

\bibitem{S-book} Z.-W. Sun, Fibonacci Numbers and Hilbert's Tenth Problem (in Chinese),
Harbin Institute of Technology Press, Harbin, 2024.

\bibitem{T85} S. P. Tung, {\it On weak number theories}, Japan. J. Math. (N.S.) {\bf 11} (1985),  203--232.

\end{thebibliography}
\end{document}